\documentclass[twoside,a4paper,12pt,centertags]{amsart}
\usepackage[T1]{fontenc}
\usepackage[utf8]{inputenc}
\usepackage{subfiles}
\makeatletter
\usepackage[section]{placeins}
\AtBeginDocument{%
	\expandafter\renewcommand\expandafter\subsection\expandafter{%
		\expandafter\@fb@secFB\subsection
	}%
}
\makeatother
\makeatletter
\@namedef{subjclassname@2020}{%
	\textup{2020} Mathematics Subject Classification}
\makeatother
\usepackage{amsfonts,amssymb,amsmath,amsthm}
\usepackage{mathrsfs}
\usepackage{nicefrac}
\usepackage{latexsym}
\usepackage{mathtools}
\usepackage[noadjust]{cite}
\usepackage{paralist}
\usepackage{longtable}
\usepackage{float}
\usepackage{pdfsync}
\usepackage{tikz, tikz-cd}
\usetikzlibrary{patterns}
\usepackage{hhline}
\usepackage{import}
\usepackage{enumitem}

\newtheorem{thm}{Theorem}[section]

\theoremstyle{definition}

\numberwithin{equation}{section}
\usepackage[plainpages=false,pdfpagelabels,
citecolor=red]{hyperref}
\usepackage{xcolor}
\hypersetup{
 colorlinks,
 linkcolor={blue!80!black},
 citecolor={red!80!black},
 urlcolor={green!40!black}
} 

\newcommand{\IC}{\mathbb{C}}
\newcommand{\IH}{\mathbb{H}}

\newcommand{\IQ}{\mathbb{Q}}
\newcommand{\R}{\mathbb{R}}

\newcommand{\cD}{\mathcal{D}} 
\newcommand{\cH}{\mathcal{H}}
\newcommand{\cI}{\mathcal{I}}

\renewcommand{\L}{\operatorname{L}} 
\newcommand{\Lloc}{\L_{\operatorname{loc}}} 
\newcommand{\C}{\operatorname{C}} 
\renewcommand{\H}{\operatorname{H}} 
\newcommand{\W}{\operatorname{W}}
\newcommand{\T}{\operatorname{T}}
\DeclareRobustCommand{\Hdot}{\dot{\H}\protect{\vphantom{H}}} 
\DeclareRobustCommand{\Wdot}{\dot{\W}\protect{\vphantom{W}}} 
\DeclareRobustCommand{\Wdot}{\dot{\W}\protect{\vphantom{W}}} 
\DeclareRobustCommand{\Lamdot}{\dot{\Lambda}\protect{\vphantom{\Lambda}}} 
\renewcommand{\S}{\mathrm{S}} 
\newcommand{\BMO}{\mathrm{BMO}}

\newcommand{\reu}{{\mathbb{R}^{1+n}_+}}

\newcommand{\tL}{\widetilde{L}}
\newcommand{\tM}{\widetilde{M}}
\newcommand{\gradx}{\nabla_x}
\renewcommand{\div}{\operatorname{div}}

\newcommand{\dnuA}{\partial_{\nu_A}} 
\newcommand{\NT}{\widetilde{N}_*} 
\newcommand{\e}{\mathrm{e}} 

\renewcommand{\d}{\, \mathrm{d}} 
\renewcommand\Re{\operatorname{Re}}

\newcommand{\Lop}{\mathcal{L}} 
\newcommand{\Le}{\mathcal{L}}
\newcommand{\cl}[1]{\overline{#1}} 
\DeclareMathOperator{\ran}{\mathsf{R}} 

\newcommand{\ind}{\mathbf{1}}
\newcommand{\dist}{\mathrm{d}}

\def\Xint#1{\mathchoice
{\XXint\displaystyle\textstyle{#1}}%
{\XXint\textstyle\scriptstyle{#1}}%
{\XXint\scriptstyle\scriptscriptstyle{#1}}%
{\XXint\scriptscriptstyle%
\scriptscriptstyle{#1}}%
\!\int}
\def\XXint#1#2#3{{\setbox0=\hbox{$#1{#2#3}{%
\int}$ }
\vcenter{\hbox{$#2#3$ }}\kern-.6\wd0}}
\def\barint{\,\Xint -} 
\def\bariint{\barint_{} \kern-.4em \barint}
\def\bariiint{\bariint_{} \kern-.4em \barint}
\renewcommand{\iint}{\int_{}\kern-.34em \int} 
\renewcommand{\iiint}{\iint_{}\kern-.34em \int} 
\title[Hardy spaces for BVP]{Hardy spaces for boundary value problems of elliptic systems with block structure}
\author{Pascal Auscher}
\address{Universit\'e Paris-Saclay, CNRS, Laboratoire de Math\'{e}matiques d'Orsay, 91405 Orsay, France \vspace{5pt}\newline\vspace{5pt}{\rm and}\newline Laboratoire Ami\'{e}nois de Math\'{e}matiques Fondamentale et Appliqu\'{e}e, UMR 7352 du CNRS, Universit\'{e} de Picardie-Jules Verne, 80039 Amiens, France}
\email{pascal.auscher@universite-paris-saclay.fr}
\author{Moritz Egert}
\address{Universit\'e Paris-Saclay, CNRS, Laboratoire de Math\'{e}matiques d'Orsay, 91405 Orsay, France}
\email{moritz.egert@universite-paris-saclay.fr}
\thanks{The authors were supported by the ANR project RAGE ANR-18-CE40-0012. This material is based upon work that got started under NSF Grant DMS-1440140 while Auscher was in residence at the MSRI in Berkeley, California, during the Spring 2017 semester. Egert also thanks this institute for hospitality.}
\date{\today}
\dedicatory{In the honor of Guido Weiss' ninetieth birthday}
\allowdisplaybreaks
\begin{document}

\begin{abstract}
We present  recent results on elliptic boundary value problems where the theory of Hardy spaces associated with operators plays a key role. 
\end{abstract}

\maketitle


\section{Introduction}

\noindent Hardy spaces have a long history in which Guido Weiss, with Elias Stein, pioneered their extension to Euclidean space based on an approach using boundary values for a Cauchy--Riemann system related to  harmonic functions \cite{steinweiss}. 

Then real analysis methods took over, first with work of Fefferman and Stein \cite{Fef-Stein},  second with the atomic theory developed by Coifman  in one dimension  \cite{coif} and Latter in any dimension \cite{L}, and third, again from the impetus of Guido Weiss, this time with Raphy Coifman, by extending the theory to spaces of homogeneous type \cite{coifw}. The presence of a special operator like the Laplacian disappears in these approaches. Still, in terms of applicability to analysis, Hardy spaces provide a continuum with the usual Lebesgue spaces $\L^p$, $1<p<\infty$, to the range $0<p\le 1$.

Numerous works treated characterizations in terms of maximal functions, atoms, molecules, square functions, Littlewood--Paley type decomposition, and even wavelet bases descriptions. Nowadays, this provides a versatile theory for treating boundedness of operators in well-understood classes, such as Calder\'on--Zygmund operators. 

However, for operators that have less regularity or different properties, these spaces may not be useful in all  their scale. For this reason, Hardy spaces associated with operators have been developed. 

One possible construction starts with vertical maximal functions for semigroups and develops an adapted atomic approach.  Let us quote the work of Dziuba\'nski and Zienkiewicz \cite{DZ1,DZ2, DZ3} for Schr\"odinger operators with certain potentials. In Auscher and Russ \cite{AR}, a link is established between Hardy spaces on domains and Hardy spaces defined via maximal functions involving  semigroups generated by elliptic operators with boundary conditions, including the Laplacian, using the atomic theory proposed by Chang et al.~\cite{CKS}.

A maximal function definition usually requires strong assumptions on the underlying operator, such as pointwise heat kernel estimates.  The next best thing is the tent space theory of Coifman et al.~\cite{CMS} based on conical square function norms. These norms incorporate $\L^2$-averages on balls 
and therefore it suffices that the operator satisfies a much weaker form of decay, called $\L^2$ off-diagonal estimates.
This machinery was developed systematically in \cite{ADMc, AMcR, HM} and inspired many subsequent works because it yields continuous scales of spaces adapted to a given operator with a well-behaved functional calculus around $\L^2$,
without requiring any $\L^p$-theory.
The theory includes atomic or molecular decompositions with adapted atoms or molecules.
Further abstractions of the atomic theory that include Hardy spaces associated with operators as a particular example have been proposed by Bernicot and Zhao \cite{BZ}.

However, one may want to see Hardy spaces associated with operators not only from an abstract point of view but investigate whether they can be identified with classical function or distribution spaces.  If this works, then it gives access to operator-valued multiplier results in limited ranges of exponents within the classical Lebesgue and Hardy spaces, because the multipliers naturally act boundedly on the adapted Hardy spaces, which are designed for that purpose.

Identifying a Hardy space associated with an operator to a classical space also yields a description of the latter in terms of a square function norm associated with the operator at hand. Such norms are useful in the treatment of boundary value problems from the harmonic analysis point of view that started from the pioneer work of Dahlberg \cite{Da}. In other words, one could go back in time and think of Hardy spaces associated with operators as a key tool box that has to appear --- implicitly or explicitly --- in the analysis of boundary value problems.  This is the leitmotif behind the program developed for elliptic boundary value problems via the first order approach  by one of us and collaborators \cite{AA1,  AusSta, AM, AA}. It proved useful in and yielded surprising developments towards boundary value problems for parabolic operators \cite{AEN1, AEN2} --- and becomes manifests in the important work of R\'osen~\cite{R, RosenJEE}, who proves  non-tangential limits and establishes the connection to classical layer potential theory. Let us stress that none of this would have been possible without the solution of Kato problem in any dimension \cite{AHLMcT} and the important Dirac operator approach \cite{AKMc}, not to mention all precursor works.

Here, we wish to announce results on well-posedness of boundary value problems for elliptic systems in block form (that is, without mixed tangential and normal derivatives) in the upper half-space with data in Lebesgue, Hardy, and homogeneous Hölder spaces. These results essentially close the topic. Prior to that, the situation was fully understood only when the boundary is one-dimensional \cite{AT}.  The admissible ranges of boundary spaces can be calculated thanks to the theory of Hardy spaces associated with a boundary operator. Details together with more exhaustive references to earlier and related literature will be presented in our forthcoming monograph \cite{AE}. 

\section{Setup}

\noindent Consider integers $n\ge 1$ and $m\ge 1$. Let $a\in \L^\infty(\R^n; \Lop(\IC^m))$ satisfy the strict accretivity condition
\begin{align}
\label{eq: accretivity a}
\Re \langle a(x) \xi, \xi \rangle \geq \lambda |\xi|^2 \quad (x \in \R^n, \, \xi \in \IC^m)
\end{align}
for some $\lambda>0$. Note that $a^{-1}$ has the same properties. Let $d\in \L^\infty(\R^n; \Lop(\IC^{nm}))$ satisfy the G\aa{}rding inequality,
\begin{align}
\label{eq: Garding d}
\Re \langle d \nabla_x v, \nabla_x v \rangle \geq \lambda \|\nabla_x v\|_2^2 \quad (v \in \C_0^\infty(\R^n; \IC^m)),
\end{align}
which in general is weaker than strict ellipticity. Here, $\langle \cdot \, , \cdot \rangle$ denotes the inner product in the respective context. One can then define a sectorial operator in $\L^2(\R^n; \IC^m)$ by $$L\coloneqq -a^{-1}\div_x d \nabla_x$$ with maximal domain $\{u\in \W^{1,2}(\R^n; \IC^m); \div_x d \nabla_x u \in \L^2(\R^n; \IC^m)\}$, where $\W^{1,2}$ denotes the usual $\L^2$-Sobolev space of order $1$. The sectoriality angle can be  any number within $[0, \pi)$, but there is a sectorial square root operator $L^{1/2}$ of angle within $[0, \nicefrac{\pi}{2})$ that generates a Poisson semigroup  $(\e^{-tL^{1/2}})_{t>0}$. This family can be seen as a solution operator for the elliptic system
\begin{equation}
\label{eq:block}
\partial_{t}(a \partial_{t}u) + \div_{x}(d\gradx  u)=0, \quad (t,x) \in (0,\infty)\times \R^n \eqqcolon \reu.
\end{equation}
Hence, the Dirichlet problem for this equation can be solved explicitly, at least for $\L^2$-data.  The regularity problem or the Neumann problem can also be solved formally. It was Kenig \cite{Ke} who observed that the interior estimates required in the classical harmonic analysis approach to these problems are linked to the Kato conjecture for $L$, that is, the nowadays known homogeneous estimate
\begin{align}
\label{eq:Kato}
\|aL^{1/2} f\|_2 \simeq \|\nabla_x f\|_2, 
\end{align}
which identifies the domain of $L^{1/2}$ as the Sobolev space $\W^{1,2}(\R^n; \IC^m)$ since $a$ is invertible in $\L^\infty(\R^n; \Lop(\IC^m))$.

Generalizing this approach to other spaces of data and obtaining existence of solutions via explicit formul\ae \, require further knowledge about the functional calculus of $L$. The exact ranges are provided by Hardy and Hardy--Sobolev spaces associated with $L$ and a related first order operator. Uniqueness of the solution, on the other hand, does not depend on Hardy space theory but on other methods relying on appropriate use of (abstract) layer potentials. Actually, it can be obtained in a larger range of exponents than the one for existence.

\section{Adapted Hardy--Sobolev spaces}

\noindent For this section, we work on $\R^n$ and functions are valued in finite dimensional complex spaces; we omit this in the notation whenever the context is clear. 

\subsection{The construction}

Given $\omega \in (0, \pi)$, we define the sector $\S_\omega^+ \coloneqq \{z\in\IC\setminus \{0\} : | \arg z|<\omega\}$ and agree on $\S_0^+\coloneqq    
(0,\infty)$. As mentioned above, $L$ is a sectorial operator in $\L^2$, which means that there is an angle $\omega \in [0, \pi)$ such that the spectrum $\sigma(L)$ is contained in $\cl{\S_\omega^+}$ (closure in $\IC$) and such that for every $\mu \in (\omega, \pi)$, we have
\begin{align}
\label{eq: constant for resolvent bound}
M_{L,\mu} \coloneqq  \sup_{z \in \IC \setminus \cl{\S_\mu^+}}  \|z(z-L)^{-1}||_{\L^2 \to \L^2} < \infty.
\end{align}
Let $\omega_L$ denote the smallest angle $\omega$ with this property.  Sectorial operators in $\L^2$ are densely defined and it can be shown that $L$ is injective and that it has dense range.
 
For injective sectorial operators, one can construct a holomorphic functional calculus to define $\psi(L)$ as a closed operator in $\L^2$ in a meaningful way if $\psi$ is holomorphic on a sector of angle larger than $\omega_L$ and has at most polynomial growth at $|z|=0$ and $|z|=\infty$. In particular, $\psi(z)=\e^{-tz^{1/2}}$ yields $\psi(L)=e^{-tL^{1/2}}$ (with the convention that we take the principal branch) and $\varphi(z)=z^{1/2}$ yields $\varphi(L)=L^{1/2}$, the square root of $L$. The particular operator $L$ has the additional property that $\psi(L)$ is bounded if $\psi$ is merely bounded. This so-called bounded $\H^\infty$-calculus is a deep result. We describe one approach in Section~\ref{subsec: Dirac} below.

The next ingredient in the construction of Hardy spaces associated with $L$ is that the resolvent family satisfies $\L^2$ off-diagonal estimates.
Let $\Omega \subseteq \IC \setminus \{0\}$. A family $(T(z))_{z \in \Omega}$ of bounded linear operators  on $\L^2$  satisfies \emph{$\L^2$ off-diagonal estimates of order $\gamma>0$} if there exists a constant $C_\gamma$ such that 
\begin{align}
\label{eq: def off-diagonal}
\|\ind_F T(z) \ind_E f \|_2 \leq C_\gamma \bigg(1+ \frac{\dist(E,F)}{|z|} \bigg)^{-\gamma} \|\ind_E f\|_2
\end{align}
for all measurable subsets $E,F \subseteq \R^n$, all $z \in \Omega$ and all $f \in \L^2$, where $\dist(E,F)$ denotes the distance between $E$ and $F$. 

The following families satisfy $\L^2$ off-diagonal estimates of arbitrarily large order:
\begin{enumerate}
	\item $((1+z^2 L)^{-1})_{z \in \S_\mu^+}$ if $\mu \in (0, \nicefrac{(\pi - \omega_L)}{2})$,
		\item $(z\nabla_x(1+z^2L)^{-1})_{z \in \S_\mu^+}$ if $\mu \in (0, \nicefrac{(\pi - \omega_L)}{2})$.
\end{enumerate}

Now, let $\psi$ be a non-zero holomorphic function on a sector $\S_\mu^+$ with $\mu > \omega_L$ that satisfies for some $\sigma,\tau>0$ the decay condition
\begin{align}
\label{eq:decay}
|\psi(z)| \lesssim (|z|^\sigma \wedge |z|^{-\tau}) \quad (z \in \S_\mu^+).
\end{align} 
Mimicking the extension to the upper half-space by convolutions in the definition of the classical Hardy spaces, one associates with $\psi$ the operator
\begin{align}
\label{eq: def Q extension}
\IQ_{\psi, L} : \cl{\ran(L)}=\L^2 \to \L^\infty(0,\infty; \L^2), \quad (\IQ_{\psi, L}f)(t) = \psi(t^2 L)f.
\end{align}
Given $s\in \R$ and $p\in (0,\infty)$, the space
\begin{align*}
\IH^{s,p}_{\psi, L} \coloneqq \{f \in \cl{\ran(L)}: \IQ_{\psi, L}f \in \T^{s,p} \}
\end{align*}
equipped with the (quasi-)norm
\begin{align*}
\|f\|_{\IH^{s,p}_{\psi,L}} \coloneqq \|\IQ_{\psi,L}f\|_{\T^{s,p}}
\end{align*}
is called \emph{pre-Hardy--Sobolev} space of smoothness $s$ and integrability $p$ adapted to $L$. The function $\psi$ is called an \emph{auxiliary function}. The space $\T^{s,p}$ is the tent space of those functions $F \in \Lloc^2(\reu)$ with finite quasi-norm
\begin{align*}
\|F\|_{\T^{s,p}} \coloneqq \|S(t^{-s} F)\|_p, 
\end{align*}
where $S$ is the conical square function
\begin{align}
\label{eq: S}
(S F)(x) \coloneqq \bigg(\iint_{|x-y|<t} |F(t,y)|^2 \, \frac{\d t \d y}{t^{1+n}} \bigg)^{\frac{1}{2}} \quad (x \in \R^n).
\end{align}

Off-diagonal estimates are used to show that all $\IH^{s,p}_{\psi, L}$ (quasi-)norms are equivalent in a certain range of parameters $\sigma,\tau$ that depends on $s, p$ and dimension.   
This provides us with a space that does not depend of the specific choice of such $\psi$. Hence, we drop $\psi$ in the notation. When $s$ and $p$ vary, they form complex interpolation scales. Moreover, holomorphic bounded functions of $L$ act continuously on such spaces. We use the suffix \emph{pre} to emphasize that they might be non complete spaces and some care is to be taken with respect to that.

\subsection{The identification problem}

For general sectorial operators, the above pre-Hardy--Sobolev spaces cannot be related to classical spaces in the expected way even when $p=2$.  But in the case of  $L$,  we have up to equivalent norms
\begin{align}
\label{eq: Hardy identification for p=2}
\IH^{0,2}_{L} = \L^2.
\end{align}
This simple equality is in fact a deep square function estimate that is, by a fundamental result due to McIntosh, equivalent to the boundedness of the $\H^\infty$-calculus.

The questions that are relevant with a view on boundary value problems for \eqref{eq:block} are 
\begin{enumerate}
	\item to identify $\IH^{s,p}_{L}$ when $s=0$ and $s=1$ and 
	\item to find the interval of boundedness for the Riesz transform $R_L\coloneqq\nabla_xL^{-1/2}$. 
\end{enumerate}
The fractional range $0<s<1$ is relevant for boundary value problems with data in fractional Sobolev spaces.

We agree that $\H^p$ denotes the classical real Hardy space on $\R^n$ when $0<p\le 1$ and the Lebesgue space $\L^p$ when $1<p<\infty$. Similarly, $\Hdot^{1,p}$ denotes the homogeneous Hardy--Sobolev space when $0<p\le 1$ and the homogeneous Sobolev space $\Wdot^{1,p}$ when $1<p<\infty$ 
with the natural quasi-norm $\|\nabla_x \cdot \|_{\H^p}$. We introduce the two sets that answer  the first question raised above:
\begin{align*}
\cH(L) &\coloneqq \big\{ p \in (1_*, \infty) : \|f\|_{\IH^{0,p}_L} \simeq \|a f\|_{\H^p} \text{ for all } f \in \L^2 \big\}, \\
\cH^1(L)&\coloneqq \big\{ p \in (1_*, \infty) : \|f\|_{\IH^{1,p}_L} \simeq \|f\|_{\Hdot^{1,p}} \text{ for all } f \in \L^2 \big\}.
\end{align*}
Here and throughout, we write $q_* \coloneqq \frac{nq}{n+q}$ for lower Sobolev conjugates. Going systematically below $p=1$ is a novelty of our approach compared to earlier references. It seems natural from the point of view of regularity theory to incorporate the possibility of having estimates in this range, as is the case for instance for operators with real coefficients when $m=1$.  
The limit exponent $1_*$ can be understood from Sobolev embeddings and duality:  the best one can hope for in absence of smoothness of the coefficients is regularity theory in H\"older spaces of exponents less than $1$.

In these regions, we can identify $\IH_L^{0,p}$ and $\IH_L^{1,p}$ with concrete quasi-normed spaces of functions, namely $$\IH_L^{0,p} = a^{-1} (\H^p \cap \L^2),$$ the image of $\H^p \cap \L^2$ under multiplication with $a^{-1}$ equipped with the image quasi-normed topology $\|\cdot\|_{a^{-1}\H^p} \coloneqq \|a \cdot\|_{\H^p}$, and  $$\IH_{L}^{1,p} = \Hdot^{1,p} \cap \L^2$$ with equivalent quasi-norm $\|\cdot\|_{\Hdot^{1,p}}$. If $p>1$, then $\|a \cdot\|_{\H^p}$ and $\|\cdot \|_{\Hdot^{1,p}}$ can be replaced with $\|\cdot\|_p$ and $\|\cdot\|_{\Wdot^{1,p}}$, respectively.
 
Likewise, the set corresponding to the second question above on the Riesz transform is
\begin{align*}
\cI(L) \coloneqq \big\{ p \in (1_*,\infty) : \|R_La^{-1}f\|_{\H^p} \lesssim \|f\|_{\H^p} \text{ for all } f \in \H^p \cap \L^2 \big\},
\end{align*}
where again the multiplication by $a^{-1}$ is only relevant when $p \leq 1$.
These three sets are intervals, and they contain $p=2$, which is the deep result mentioned above. Using extended Calder\'on--Zygmund theory together with functional calculus, 
one can calculate in a sharp fashion $\cH(L), \cI(L)$ and the upper endpoint of $\cH^1(L)$ as functions of \emph{four critical numbers for $L$}:

 \begin{thm}
\label{thm: characterization} One has
\begin{enumerate}
\item $ \cH(L) = (p_-(L), p_+(L)),$
\item $ \cH^1(L) \supseteq (q_-(L)_* \vee 1_*, q_+(L)),$
\item $ \cI(L) = (q_{-}(L), q_{+}(L)).$
\end{enumerate}
The upper bound in (ii) is sharp and not attained. Moreover, $aL^{1/2}$ extends to a bounded operator $\Hdot^{1,p} \to \H^p$ if $p\in \cH(L) \cup \cH^1(L) $ and is an isomorphism if and only if $p\in \cI(L)$.
\end{thm}

\subsection{The critical numbers}

The definition of these four numbers is given in terms of bounds for the resolvents of $L$ and their gradients: 

\begin{itemize}
	\item  $(p_{-}(L), p_{+}(L))$ is the maximal open set within  $(1_*, \infty)$ for which the family $(a(1+t^2L)^{-1}a^{-1})_{t>0}$ is uniformly bounded on $\H^p$. \\[-8pt]
	\item  $(q_{-}(L),q_{+}(L))$ is the maximal open set within $(1_*, \infty)$ for which the family $(t\nabla_{x}(1+t^2L)^{-1}a^{-1})_{t>0}$ is uniformly bounded on $\H^p$. 
\end{itemize}

\smallskip

It is remarkable that simple boundedness information on the resolvent family completely rules the Hardy space theory of $L$ and the boundedness of its Riesz transform.  It can be shown that replacing the resolvent by the Poisson semigroup $(\e^{-tL^{1/2}})_{t>0}$ in the definition above leads to the same critical numbers. If the sectoriality angle of $L$ exceeds $\nicefrac{\pi}{2}$, then one has no "heat" semigroup and this is why we have to deviate from most of the literature on the topic right from the start.

The critical numbers characterize many other boundedness properties related to $L$ and in fact there are only three of them since one can show
$q_-(L)=p_-(L)$. Moreover, essentially due to Sobolev embeddings, one has $p_+(L)\ge q_+(L)^*$, where $p^*=\frac{np}{n-p}$ if $p<n$ and $p^* = \infty$ otherwise. The best general conclusion for the critical numbers is
\begin{align*}
(p_{-}(L), p_{+}(L)) 
&\supseteq 
{\begin{cases}
(\frac{1}{2}, \infty) & \text{if $n=1$} \\
[1, \infty) & \text{if $n = 2$} \\
[\frac{2n}{n+2}, \frac{2n}{n-2}] & \text{if $n \geq 3$}
\end{cases}}
\intertext{and}
(q_{-}(L), q_{+}(L)) &\supseteq
{\begin{cases}
(\frac{1}{2}, \infty) & \text{if $n=1$} \\
[\tfrac{2n}{n+2}, 2] & \text{if $n \ge 2$} \\
\end{cases}}\ .
\end{align*}
In specific cases, one  can say more. If $L$ has constant coefficients, then $q_-(L)=p_-(L)=1_*$ and 
$q_+(L)=p_+(L)=\infty$. We have the same numbers for $L=-a^{-1} \Delta_x$ (that is, $d$ is the identity) or when $n=1$ with arbitrary $d$. If $n= 2$ or if $n\ge3$ and $d$ is  a real valued $(n \times n)$-matrix (hence, $m=1$), then $p_+(L)=\infty$ and 
$q_-(L)=p_-(L)<1$. 

\subsection{The underlying perturbed Dirac operator}
\label{subsec: Dirac}

The interior control for the boundary value problems requires us to look at expressions such as $\nabla_x \e^{-tL^{1/2}}f$, which means that we have to leave the functional calculus of $L$. Such expressions can often be handled in a direct way by specific arguments but they can also be understood in a unified way within the functional calculus of a perturbed first order Dirac operator $DB$, where 
\begin{align*}
D\coloneqq\begin{bmatrix} 0 & \div_x  \\ 
-\nabla_x & 0 \end{bmatrix},
\quad
B\coloneqq\begin{bmatrix} a^{-1} & 0  \\ 
	0 & d \end{bmatrix}.
\end{align*} 
This operator $DB$ is not sectorial  but bisectorial in $\L^2(\R^n; \IC^{m} \times \IC^{mn})$, which means that its spectrum is contained in a double sector $\cl{\S_{\omega_{DB}}^+} \cup -\cl{\S_{\omega_{DB}}^+}$, with $\omega_{DB}<\nicefrac{\pi}{2}$  and one has resolvent estimates of the form \eqref{eq: constant for resolvent bound} away from it.  In particular 
\begin{align}
\label{eq: tL and tM}
(DB)^2=\begin{bmatrix} -\div_x d \nabla_x a^{-1 \vphantom{\tL}}& 0 \\ 0 & -\nabla_x a^{-1} \div_x d \vphantom{\tM}\end{bmatrix}
\end{align} 
is a sectorial operator in $\L^2$ with angle not exceeding $2\omega_{DB}$ that "contains" an operator similar to $L$ under conjugation with $a$. The square root $[DB]\coloneqq ((DB)^2)^{1/2}$ generates a semigroup that contains in the same manner a family similar to the Poisson semigroup for $L$. Then, for example,
\begin{equation}\label{eq:link}
DB\e^{-t[DB]}\begin{bmatrix} af   \\ 
 0 \end{bmatrix}= - \begin{bmatrix} 0   \\ 
 \nabla_x \e^{-tL^{1/2}}f \end{bmatrix} 
\end{equation}
is a natural object in the calculus for $DB$ and one has an interest in introducing adapted Hardy--Sobolev spaces $\IH^{s,p}_{\psi,DB}$. 

Indeed,  the family $((1+itDB)^{-1})_{ t\in \R\setminus \{0\}}$ satisfies $\L^2$ off-diagonal estimates and this allows us to develop the theory similarly to the one for $L$ with two main differences. First, as $DB$ is not injective (except in dimension $n=1$),  we have to work systematically on the $\L^2$-closure of its range, which is in topological direct sum with its nullspace. Its range is the same as the one for $D$, denoted $\ran(D)$; this independence with respect to $B$ is a useful fact.  Secondly, we use first-order scaling in $t$ for the extension operators $\IQ_{\psi,DB}$ to make sure that $s$ has the same meaning as smoothness parameter as for the $L$-adapted spaces. When $s=0$, let us write $\IH^{p}_{DB}$ for $\IH^{0,p}_{\psi,DB}$.

The equality 
\begin{align}
\IH^{2}_{DB} = \cl{\ran(D)}
\end{align}
is again a deep square function estimate that is equivalent (by McIntosh's result) to the boundedness of the $\H^\infty$-calculus for $DB$ on $\L^2$. The latter implies \eqref{eq:Kato} and the boundedness of the $\H^\infty$-calculus for $(DB)^2$ on $\L^2$, which due to the block structure is equivalent to the $\H^\infty$-calculus for $L$ on both $\L^2$ and $\Wdot^{1,2}$. Thus, it is significantly more information than \eqref{eq: Hardy identification for p=2}.

For $B=I$, one can identify $
\IH^{p}_{D}$ with $\H^p \cap\, \cl{\ran(D)}$ when $p\in (1_*, \infty)$, using Fourier multipliers. Hence, there is a natural identification region, which contains $p=2$, 
$$
\cH(DB) \coloneqq \big\{ p \in (1_*, \infty) : \|h\|_{\IH^p_{DB}} \simeq \|h\|_{\H^p} \text{ for all } h \in \cl{\ran(D)}  \big\}.
$$
Heuristically, the theory for $\IH_{DB}^p$ comprises the theory for $L$ at both smoothness scales $s=0$ and $s=1$. On the level of Hardy spaces, this can be expressed through the identity $\cH(DB) = \cH(L) \cap \cH^1(L)$, so that $$\cH(DB)= (q_{-}(L),q_{+}(L)).$$ 
On the level of boundary value problems, it leads us to the principle that the $DB$-adapted theory applies to  all boundary value problems of Neumann and Dirichlet type, whereas the $L$-adapted theory allows one to separate issues in the sense that $s=0$ corresponds to the Dirichlet problem and $s=1$ corresponds to the regularity problem.
\section{The Dirichlet  problem}

\noindent  For harmonic functions it is well known from the works of Calder\'on~\cite{C} and Stein~\cite{Stein70} that non-tangential control  gives access to almost everywhere convergence at the boundary.  In the context of elliptic equations with real and measurable coefficients,
Dahlberg \cite{Da} then formulated the Dirichlet problem with  non-tangential maximal estimates in order to recover almost everywhere boundary limits. 

Since  weak solutions  for general systems might not be regular, we use the Whitney average variants of the non-tangential maximal function in order to pose our boundary value problems. When getting back to systems where solutions have pointwise values, these variants turn out to be equivalent to the usual pointwise control. More precisely, we let
\begin{align*}
\NT(F)(x):= \sup_{t>0}  \bigg(\bariint_{W(t,x)} |F(s,y)|^2\, {\d s \d y}\bigg)^{1/2} \quad (x\in \R^n),
\end{align*}
where $W(t,x) \coloneqq (\nicefrac{t}{2}, 2t) \times B(x,t)$ is a Whitney region.

For  $1<p<\infty$, the  Dirichlet problem with non-tangential maximal control  and data $f\in \L^p(\R^{n}; \IC^m) $ consists in solving 
\begin{equation*}
(D)_{p}^\Le  \quad\quad
\begin{cases}
\partial_{t}(a \partial_{t}u) + \div_{x}(d\gradx  u)=0   & (\text{in }\reu), \\
\NT (u)\in \L^p(\R^n),   \\
\lim_{t \to 0} \bariint_{W(t,x)} |u(s,y)-f(x)| \, \d s \d y = 0 & (\text{a.e. } x\in \R^n).\end{cases}
\end{equation*}
The problem for $p=1$ is formulated analogously with $a^{-1}\H^1$ as data space on the boundary. We note that $\H^1 \subseteq \L^1$ and therefore also the image $a^{-1}\H^1$ is a Banach space that embeds into $\L^1$.

As usual, well-posedness means existence, uniqueness, and continuous dependence on the data. Compatible well-posedness means well-posedness together with the fact that the solution agrees with the energy solution when the data $f$ also belong to the trace space $\Hdot^{1/2}(\R^n;\IC^m)$.  The energy solution is  constructed  in the homogeneous Sobolev space $\Wdot^{1,2}(\reu; \IC^m)$ via the Lax--Milgram lemma.

The following result, whose range of exponents is likely to be optimal, summarizes the situation.  

\begin{thm}
\label{thm: blockdir} 
Let  $p \geq 1$ be such that $p_{-}(L)< p <  p_+(L)^*$. Given $f\in \L^p(\R^n; \IC^m)$ when $p>1$ and $f\in a^{-1} \H^1(\R^n; \IC^m)$ when $p=1$, the Dirichlet problem $(D)_{p}^{\Le}$ is compatibly well-posed. The solution $u$ has the following additional properties.  
\begin{enumerate}
	\item There is comparability
	\begin{align*} 
	\|\NT(u)\|_{p} \simeq \|af\|_{\H^p} \simeq \|S (t\nabla u)\|_{p}.
	\end{align*}
	\item The non-tangential convergence improves to $\L^2$-averages
	\begin{align*}
	\qquad \lim_{t \to 0} \bariint_{W(t,x)} |u(s,y)-f(x)|^2 \, \d s \d y = 0  \quad (\text{a.e. } x\in \R^n).
	\end{align*}
	\item When $p<p_+(L)$, then   $au$ is of class\footnote{As usual, the notation $\C_0([0,\infty))$ means continuity and limit $0$ at infinity.} $\C_0([0,\infty); \H^p(\R^n; \IC^m)) \cap \C^\infty((0,\infty); \H^p(\R^n; \IC^m))$ with $au(0, \cdot) = af$ and 
	\begin{align*}
	\sup_{t>0} \|au(t,\cdot)\|_{\H^p} \simeq \|af\|_{\H^p}.
		\end{align*}
	\item When $p \geq p_+(L)$, then for all $T>0$ and compact $K \subseteq \R^n$, $u$ is of class $\C([0,T]; \L^2(K; \IC^m))$ with $u(0, \cdot) = f$ and there is a constant $c=c(T,K)$ such that
	\begin{align*}
	\sup_{0<t \leq T} \|u(t,\cdot)\|_{\L^2(K)} \lesssim c \|f\|_p.
	\end{align*}	
\end{enumerate}
\end{thm}

As expected, the above solution is provided by $u(t,x) =  \e^{-tL^{1/2}} f(x)$ for $f \in a^{-1}(\H^{p} \cap \L^2)$ and  by an extension by density of this expression in  the respective topologies for  $f \in a^{-1}\H^{p}$. In the range $p<p_+(L)$ one can use the extension to a proper $\C_0$-semigroup on adapted Hardy spaces, which explains the regularity result (iii) given the identification range $\cH(L)$ in Theorem~\ref{thm: characterization}.  However, and this was never observed before, the range of exponents in the statement exceeds by one Sobolev exponent the range provided by the semigroup theory. This means that in this case, $u$ is extended as a function of both variables $t$ and $x$ simultaneously.

The link to adapted Hardy spaces is seen in the square function estimates in (i).  For $p>2$, some specific considerations relying on off-diagonal estimates are used together with a generic result in $\L^p$-functional calculus from \cite{CDMcY}. 
The argument for $p<2$ is more involved and relies on adapted Hardy spaces. Indeed, note that $t\nabla u$ has two components. The first one is $t\partial_t u$, so for $u(t,\cdot)= \e^{-tL^{1/2}} f$ one must check that $\psi(z)= z^{1/2}\e^{-z^{1/2}}$ is an admissible  auxiliary function for $\IH^{0,p}_{L}$. The other one is $t\nabla_x u$, in which case one can use \eqref{eq:link} and the $DB$-theory.

Finally, the non-tangential maximal function estimates in (i) can be obtained from the general theory of $L$-adapted Hardy spaces when $p\le 2$. The case $p>2$ is more straightforward as the almost everywhere limit in (ii) furnishes the lower bound. 

From the point of view of $L$-adapted Hardy spaces, the natural boundary space for the Dirichlet problem is the completion of $\IH_L^p = a^{-1}(\H^p \cap \L^2)$ for the quasi-norm $\|a \cdot \|_{\H^p}$. Since $a$ is not smooth, there is no reason to believe that a completion can be realized within the ambient space of distributions if $p<1$ and hence there is no meaningful way of posing a Dirichlet problem that is compatible with the Lax--Milgram approach. What remains true though, again thanks to the theory of adapted Hardy spaces, are the estimates in the theorem above when $u$ is the Poisson semigroup extension of data in $a^{-1}(\H^p \cap \L^2)$. 

\begin{thm}
\label{thm: blockdirHp} 
Suppose $p_{-}(L)< 1$ and let $p_-(L)<  p < 1$. Given $f\in a^{-1}(\H^p(\R^n; \IC^m) \cap \L^2(\R^n; \IC^m))$, the semigroup extension $u(t,x) =  \e^{-tL^{1/2}} f(x)$ satisfies (i), (ii), (iii) of Theorem~\ref{thm: blockdir}.
\end{thm}

\section{The Dirichlet  problem for H\"older data}

\noindent One may wonder whether in the case $p_+(L) > n$ one can extend the results for the Dirichlet problem to exponents "beyond $\infty$", which we think of corresponding to the homogeneous Hölder spaces $\Lamdot^\alpha(\R^n;\IC^m)$, $0 \leq \alpha < 1$, with the endpoint case $\Lamdot^0 \coloneqq \BMO$. In this case, we choose  the interior control as the Carleson functional
\begin{align*}
C_\alpha F(x) \coloneqq \sup_{t>0} \frac{1}{t^\alpha} \bigg(\frac{1}{t^n} \int_0^{t} \int_{B(x,t)} |F(s,y)|^2 \frac{\d y \d s}{s} \bigg)^{1/2}.
\end{align*}
For $\alpha\in [0,1)$, the Dirichlet problem with data  $f\in \Lamdot^\alpha(\R^n;\IC^m)$ consists in solving 
\begin{equation*}
(D)_{\Lamdot^\alpha}^\Le  \quad\quad
\begin{cases}
\partial_{t}(a \partial_{t}u) + \div_{x}(d\gradx  u)=0   & (\text{in } \reu), \\
C_\alpha (t\nabla u)\in \L^\infty(\R^n),   \\
\lim_{t \to 0} \barint_{t}^{2t} |u(s,\cdot)-f| \, \d s = 0 & (\text{in } \Lloc^2(\R^n; \IC^m)).\end{cases}
\end{equation*}
The choice of Carleson control for the problem with $\BMO$ data is completely natural in view of the Fefferman--Stein characterization of $\BMO$ by the Carleson functional $C_0$ for harmonic functions. Making this choice for $\Lamdot^\alpha$ when $\alpha>0$ is a convenience in terms of unification of the presentation. Other choices, which are weaker in terms of interior control, could be taken;  
 we could  use Whitney average convergence at the boundary and uniform control of Whitney averages for $u-f$ on $W(t,x)$ by $t^\alpha$  and  all results below remain the same. On the other hand, the $\BMO$ Dirichlet problem  with Carleson control and Whitney average convergence at the boundary appears to be out of reach.

For the next  statement,  $L^\sharp$ denotes the boundary operator for the adjoint equation $\partial_{t}(a^* \partial_{t}u) + \div_{x}(d^*\gradx  u)= 0$, that is,
 \begin{equation*}
L^\sharp= -(a^*)^{-1} \div_x d^* \nabla_x.
\end{equation*}

\begin{thm}
\label{thm: Holder-dir}
 Suppose that $p_+(L) > n$ and that $0 \leq \alpha < 1 - \nicefrac{n}{p_+(L)}$. Then the Dirichlet problem $(D)_{\Lamdot^\alpha}^{\Le}$ is compatibly well-posed. Given $f\in \Lamdot^\alpha(\R^n; \IC^m)$, the unique solution $u$ has the following additional properties.
\begin{enumerate}
	\item There is comparability
	\begin{align*}
	\|C_\alpha(t \nabla u)\|_\infty \simeq \|f\|_{\Lamdot^\alpha}.
	\end{align*}
	\item   $u$ is of class $\C([0,T]; \L^2(K; \IC^m))$ with $u(0, \cdot) = f$ for every $T>0$ and compact $K \subseteq \R^n$.
	\item If moreover $p_{-}(L^\sharp) < 1$ and $\alpha < n(\nicefrac{1}{p_-(L^\sharp)} -1)$, then $u$ is of class $\C_0([0,\infty); \Lamdot^\alpha_{\text{weak}^*}(\R^n; \IC^m)) \cap \C^\infty((0,\infty); \Lamdot^\alpha_{\text{weak}^*}(\R^n; \IC^m))$ and
	\begin{align*}
	\sup_{t>0} \|u(t,\cdot)\|_{\Lamdot^\alpha} \simeq \|f\|_{\Lamdot^\alpha}.
	\end{align*}
	In addition, $u$ is of class $\Lamdot^\alpha(\overline{\reu};\IC^m)$, with 
	\begin{align*}
	\|u\|_{\Lamdot^\alpha(\overline{\reu})}\lesssim \|f\|_{\Lamdot^\alpha}.
	\end{align*}
\end{enumerate}
\end{thm}

Since $\Lamdot^\alpha \cap \L^2$ is \emph{not} dense in $\Lamdot^\alpha$, we cannot extend the Poisson semigroup to the boundary space by density. In (iii), $\Lamdot^\alpha$ is considered as the dual space of $\H^p$, where $\alpha = n(\nicefrac{1}{p} -1)$, with the weak$^*$ topology. The assumption in (iii) implies $p_+(L)=\infty$ and that the solution can be constructed by duality, using the extension of the Poisson semigroup for $L^* = a^*L^\sharp(a^*)^{-1}$ to $\H^p$. Therefore, the solution keeps the $\Lamdot^\alpha$-regularity in the interior. This construction has appeared earlier. 

The construction of the solution under the mere assumption that $p_+(L) > n$ is much more general and we have
\begin{align*}
	u(t,x) = \lim_{j \to \infty} \e^{-t L^{1/2}} (\ind_{\{|\,\cdot\,| < 2^j\}} f)(x),
\end{align*}
where we have to use the assumption $p_+(L) > n$ already to prove convergence of the right-hand side in $\Lloc^2(\reu; \IC^m)$. This opens the possibility of uniquely solving Dirichlet problems for H\"older continuous (or $\BMO$) data, while producing solutions that have no reason to be in the same class in the interior of the domain. To the best of our knowledge, this phenomenon is observed for the first time. Note also that $p_+(L) > n$ always holds in dimension $n \leq4$, so that in these dimensions, the $\BMO$ Dirichlet problem is compatibly well-posed.
\section{The Regularity problem}

\noindent For $1_*<p<\infty$ the  regularity problem with data  $f\in \Hdot^{1,p}(\R^n;\IC^m)$ consists in solving 
\begin{equation*}
(R)_{p}^\Le  \quad\quad
\begin{cases}
\partial_{t}(a \partial_{t}u) + \div_{x}(d\gradx  u)=0   & (\text{in } \reu), \\
\NT(\nabla u)\in \L^p(\R^n),    \\
\lim_{t \to 0} \bariint_{W(t,x)} |u(s,y)-f(x)| \, \d s \d y = 0 & (\text{a.e. } x\in \R^n).\end{cases}
\end{equation*}
By Hardy--Sobolev embeddings, we have $f \in \Lloc^1$ so that the convergence at the boundary is meaningful. 

Our result is as follows, with a range of exponents that corresponds to $\cH^1(L)$.

\begin{thm}
\label{thm: blockreg} 
Let  $(q_{-}(L)_*\vee 1_*) < p < q_+(L)$. The regularity problem $(R)_{p}^{\Le}$ is compatibly well-posed. Given $f\in \Hdot^{1,p}(\R^n; \IC^m)$, the unique solution $u$ has the following additional properties.
\begin{enumerate}
	\item There are estimates
	\begin{align*}
	\hspace{1.5cm} \|\NT(\nabla u)\|_{p} \simeq  \|S (t\nabla \partial_{t}u)\|_{p} \simeq  \|\nabla_x f\|_{\H^p} \gtrsim \|g\|_{\H^p}
	\end{align*}
	with $g \coloneqq -a{L^{1/2}}f$ being the conormal derivative of $u$, where the square root extends from $\Hdot^{1,p}(\R^n; \IC^m) \cap \W^{1,2}(\R^n; \IC^m)$ by density. 
	\item For a.e.\ $x \in \R^n$ and all $t>0$,
	\begin{align*}
		\bigg(\bariint_{W(t,x)} |u(s,y) - f(x)|^2 \, \d s \d y \bigg)^{\frac{1}{2}} \lesssim t \NT(\nabla u)(x).
	\end{align*}
	In particular, the non-tangential convergence improves to $\L^2$-averages. {Moreover, $\lim_{t \to 0} u(t,\cdot) = f$ in $\cD'(\R^n)$.}
\item If $p \geq 1$, then for a.e.\ $x \in \R^n$,
	\begin{align*}
	\lim_{t \to 0} \bariint_{W(t,x)} \bigg|\begin{bmatrix} a\partial_t u \\ \nabla_x u \end{bmatrix} - \begin{bmatrix} g(x) \\ \nabla_x f(x) \end{bmatrix}\bigg|^2 \, \d s \d y = 0,
	\end{align*}
	where $g$  is as in (i).
	 	\item $u$ is of class $\C_{0}([0,\infty); \Hdot^{1,p}(\R^n; \IC^m)) \cap \C^\infty((0,\infty); \Hdot^{1,p}(\R^n; \IC^m))$ with 
	$$\|\nabla_x f\|_{\H^p} \simeq \sup_{t > 0} \|\nabla_{x}u(t, \cdot )\|_{\H^p}.$$
	If $p<n$, then up to a constant\footnote{The constant is chosen via Hardy--Sobolev embeddings such that $f\in\L^{p^*}$.}  $u \in \C_{0}([0,\infty); \L^{p^*}(\R^n; \IC^m)) \cap \C^\infty((0,\infty);  \L^{p^*}(\R^n; \IC^m))$ with $u(0,\cdot) = f$ and
	\begin{align*}
	\|f\|_{p^*} \leq \sup_{t > 0} \|u(t, \cdot )\|_{p^*} \lesssim \|\nabla_x f\|_{\H^p} + \|f\|_{p^*}.
	\end{align*}
	\item If $p>p_{-}(L)$, then $a\partial_{t}u$ is of class $ \C_{0}([0,\infty); \H^p(\R^n; \IC^m))$ and, with $g$ as in (i),
	\begin{align*}
	\qquad \|\NT(\partial_{t} u)\|_{p} \simeq \sup_{t\ge 0} \|a\partial_{t}u(t, \cdot )\|_{\H^p} \simeq \|g\|_{\H^p}\simeq \|\nabla_x f\|_{\H^p}.
	\end{align*}
\end{enumerate}
\end{thm}

As we have $p_-(L)=q_-(L)$, the range of exponents exceeds the semigroup range for $L$ by up to one Sobolev exponent below. 

Again, the bulk of the result is thanks to the Hardy--Sobolev theory, this time for $\IH^{1,p}_L$, and also due to the observation that the spatial gradient intertwines  $\IH^{1,p}_L$ and $\IH^{0,p}_{\widetilde M}$, with $\widetilde M= -  \nabla_x a^{-1}\div_x d$, via the formal relation $\nabla_xL=\widetilde M\nabla_x$. In fact, $\widetilde M$ is the lower right corner of $(DB)^2$ in \eqref{eq: tL and tM}. Likewise, the non-tangential maximal estimates and limits rely on the link between the Poisson semigroup for  $L$ and the semigroup for $[DB]$ via the intertwining property
\begin{align*}
\e^{-t [DB]} \begin{bmatrix} -a{L^{1/2}}f \\ \nabla_x f \end{bmatrix} 
= \begin{bmatrix} a\partial_t  \e^{-t {L^{1/2}}}f  \\ \nabla_x \e^{-t {L^{1/2}}}f \end{bmatrix}.
\end{align*}
\section{The Neumann problem}

\noindent For $1_* < p < \infty$, the  Neumann problem with data $g\in \H^p(\R^n; \IC^m)$ consists in solving (modulo constants)
\begin{equation*}
(N)_{p}^\Le  \quad\quad
\begin{cases}
\partial_{t}(a \partial_{t}u) + \div_{x}(d\gradx  u)=0, \quad &  (\text{in }\reu), \\
\NT (\nabla u)\in \L^p(\R^n),   \\
\lim_{t \to 0} a \partial_t u(t, \cdot) = g & (\text{in } \cD'(\R^n; \IC^m)).\end{cases}
\end{equation*}
Note that due to the block structure, $a \partial_t u$ is indeed the (inward) conormal derivative $\dnuA u = e_0 \cdot A \nabla u$. 

Here, the range of exponents $p$ is the set $\cH(DB)$ for the identification of $\IH_{DB}^p$ and in our block case, this is same as the range of boundedness of the Riesz transform.

\begin{thm}
\label{thm: blockneu}
Let  $q_-(L) < p < q_+(L)$. Then the Neumann problem $(N)_p^{\Le}$ is compatibly well-posed (modulo constants). Given a Neumann data $g\in \H^p(\R^n; \IC^m)$ and defining $f \coloneqq -(aL^{1/2})^{-1}g \in \Hdot^{1,p}(\R^n)$, all the properties for the solution listed in  Theorem \ref{thm: blockreg} hold.
\end{thm}

In the case $p\ge 1$ one can pose $(N)_p^{\Le}$ with convergence of Whitney averages 
$$\lim_{t \to 0} \bariint_{W(t,x)} a \partial_tu \, \d s \d y = g(x)  \qquad (\text{a.e. } x\in \R^n)
$$
instead of convergence in the sense of distributions at the boundary 
and still obtain a compatibly well-posed problem with the same solution as before.

We use the $\H^p$-boundedness of $\nabla_x (aL^{1/2})^{-1} = R_L a^{-1}$ to convert the Neumann data $g$ into an appropriate regularity data $f$. In other words, we use the ansatz $u = - \e^{-t L^{1/2}} (aL^{1/2})^{-1} g = \e^{-t L^{1/2}} f$. This is why we need to stay withing the interval $\cI(L)$ for the Riesz transform and obtain immediately the same estimates as in Theorem \ref{thm: blockreg}.

Compatible well-posedness as above has been obtained in \cite{AM} for elliptic systems, not necessarily in block form, in a certain range $I_L$ of exponents where $\IH^{p}_{DB}=\IH^p_D$ and certain additional technical conditions hold. The strategy in that paper begins by writing \eqref{eq:block} in the equivalent form
\begin{align*}
\partial_t \begin{bmatrix} a \partial_t u \\ \nabla_x u \end{bmatrix}
+ \begin{bmatrix} 0& \div_x \\ -\nabla_x & 0 \end{bmatrix} 
\begin{bmatrix} a^{-1} & 0  \\ 0 & d \end{bmatrix}
\begin{bmatrix} a \partial_t u \\ \nabla_x u \end{bmatrix}
= \begin{bmatrix}0 \\ 0 \end{bmatrix},
\end{align*}
where the second line is a dummy equation and where $F \coloneqq [a \partial_t u, \nabla_x u]^\top$ is called \emph{conormal gradient}. Hence, $\partial_t F + DB F = 0$. The main thesis there is that if $p \in I_L$, then any solution $u$ with $\NT (\nabla u)\in \L^p(\R^n)$ is obtained by the evolution of the $[DB]$-semigroup on a vector $F_0$ in the positive spectral space for $DB$ within the completion of  $\IH^{p}_{DB}=\IH^p_D$ in $\H^p$.  At the boundary we should have $F|_{t=0} = [g, f]^\top$, which yields for each $t>0$ the representation
\begin{align*}
\begin{bmatrix} a \partial_t u \\ \nabla_x u \end{bmatrix}=\e^{-t[DB]}\begin{bmatrix} g \\ f \end{bmatrix}  
=  \begin{bmatrix} aL^{1/2} \e^{-t L^{1/2}} (aL^{1/2})^{-1} g  \\ -\nabla_x\e^{-t L^{1/2}} (aL^{1/2})^{-1} g \end{bmatrix}.
\end{align*}
From this,  we see that $u$ has to be as in the ansatz.

In view of these results, our new contribution here purely is a theorem on Hardy spaces: in the block case we identify the range $I_L$ in which all of the above works as the full interval $\cH(DB) = (q_-(L), q_+(L))$.


\end{document}